\documentclass[twosided,reqno]{amsart}

\usepackage{amsmath}
\usepackage{amssymb}
\usepackage{amsthm}

\theoremstyle{theorem}
\newtheorem{theorem}{\scshape Theorem }[section]

\newtheorem{corollary}[theorem]{\scshape Corollary}
\newtheorem{proposition}[theorem]{\scshape Proposition}
\theoremstyle{definition}

\numberwithin{equation}{section}

\begin{document}

\title[Sheffer sequences of polynomials and their applications]{Sheffer sequences of polynomials and their applications}

\author{Taekyun Kim$^1$}
\address{$^1$ Department of Mathematics, Kwangwoon University, Seoul 139-701, Republic of Korea.}
\email{tkkim@kw.ac.kr}

\author{Dae San Kim$^2$}
\address{$^2$ Department of Mathematics, Sogang University, Seoul 121-742, Republic of Korea.}
\email{dskim@sogang.ac.kr}

\author{Seog-Hoon Rim$^3$}
\address{$^3$ Department of Mathematics Education, Kyungpook National University, Taegu 702-701, Republic of Korea.}
\email{shrim@knu.ac.kr}

\author{Dmitry V. Dolgy$^4$}
\address{$^4$ Hanrimwon, Kwangwoon University, Seoul 139-7011, Republic of Korea.}
\email{dgekw2011@gmail.net}


\maketitle

\begin{abstract}
In this paper, we investigate some properties of several Sheffer sequences of several polynomials arising from umbral calculus. From our investigation, we can derive many interesting identities of several polynomials.
\end{abstract}

\section{Introduction}

As is well known, the {\it{Bernoulli polynomials}} of order $a$ are defined by the generating function to be
\begin{equation}\label{1}
\left(\frac{t}{e^t-1}\right)^a e^{xt}=\sum_{n=0} ^{\infty}B_n ^{(a)}(x)\frac{t^n}{n!},{\text{ (see [1-10])}},
\end{equation}
and the {\it{Narumi polynomials}} are also given by
\begin{equation}\label{2}
\left(\frac{\log (1+t)}{t}\right)^a(1+t)^x=\sum_{n=0} ^{\infty}\frac{N_n ^{(a)}(x)}{n!}t^n,{\text{ (see [18,19])}}.
\end{equation}
In the special case, $x=0$, $N_n ^{(a)} (0)=N_n ^{(a)}$ are called the {\it{Narurni numbers}}.

Throughout this paper, we assume that $\lambda\in{\mathbb{C}}$ with $\lambda \neq 1$. {\it{Frobenius-Euler polynomials}} of order $a$ are defined by the generating function to be
\begin{equation}\label{3}
\left(\frac{1-\lambda}{e^t-\lambda}\right)^a e^{xt}=\sum_{n=0} ^{\infty} H_n ^{(a)}(x|\lambda)\frac{t^n}{n!},{\text{ (see [10-21])}}.
\end{equation}

The {\it{Stirling number of the second kind}} is also defined by the generating function to be
\begin{equation}\label{4}
(e^t-1)^n=n!\sum_{k=n} ^{\infty}S_2(k,n)\frac{t^k}{k!},{\text{ (see [9,10,18,19])}},
\end{equation}
and the {\it{Stirling number of the first kind}} is given by
\begin{equation}\label{5}
(x)_n=x(x-1)\cdots(x-n+1)=\sum_{l=0} ^n S_1(n,l)x^l, {\text{ (see [9,11,18,19])}}.
\end{equation}

Let
\begin{equation}\label{6}
{\mathcal{F}}=\left\{ \left.f(t)=\sum_{k=0} ^{\infty} \frac{a_k}{k!} t^k~\right|~ a_k \in {\mathbb{C}} \right\}.
\end{equation}
Let ${\mathbb{P}}$ be the algebra of polynomials in the variable $x$ over ${\mathbb{C}}$ and ${\mathbb{P}}^{*}$ be the vector space of all linear functionals on ${\mathbb{P}}$. The action of the linear functional $L$ on a polynomial $p(x)$ is denoted by $\left<L|p(x)\right>$. We recall that the vector space structures on ${\mathbb{P}}^{*}$ are defined by $\left< L+M|p(x)\right>=\left<L|p(x)\right>+\left<M|p(x)\right>$, $\left< cL|p(x)\right>=c\left<L|p(x)\right>$, where $c$ is a complex constant (see \cite{18,19}).

For $f(t)=\sum_{k=0} ^{\infty} a_k\frac{t^k}{k!} \in {\mathcal{F}}$, we define a linear functional $f(t)$ on ${\mathbb{P}}$ by setting
\begin{equation}\label{7}
\left<f(t)|x^n \right>=a_n,~(n\geq 0).
\end{equation}
By \eqref{6} and \eqref{7}, we get
\begin{equation}\label{8}
\left<t^k | x^n \right>=n! \delta_{n,k},~(n,k \geq 0),
\end{equation}
where $\delta_{n,k}$ is the Kronecker symbol (see \cite{09,10,11,18,19}).

Suppose that $f_L(t)=\sum_{k=0} ^{\infty} \frac{\left<L|x^k\right>}{k!}t^k$. Then we have $\left<f_L(t)|x^n\right>=\left<L|x^n\right>$ and $f_L(t)=L$.
Thus, we note that the map $L \mapsto f_L(t)$ is a vector space isomorphism from ${\mathbb{P}}^{*}$ onto ${\mathcal{F}}$. Henceforth, ${\mathcal{F}}$ will be thought of as both a formal power series and a linear functional. We shall call ${\mathcal{F}}$ the {\it{umbral algebra}}. The umbral calculus is the study of umbral algebra (see \cite{09,10,11,18,19}).

The order $o(f(t))$ of the non-zero power series $f(t)$ is the smallest integer $k$ for which the coefficient of $t^k$ does not vanish. If $o(f(t))=1$, then $f(t)$ is called a {\it{delta series}}. If $o(f(t))=0$, then $f(t)$ is called an {\it{invertible series}}. Let $o(f(t))=1$ and $o(g(t))=0$. Then there exists a unique sequence $S_n(x)$ of polynomials such that $\left<g(t)f(t)^k|S_n(x)\right>=n!\delta_{n,k}$ $(n,k \geq 0)$. The sequence $S_n(x)$ is called {\it{Sheffer sequence}} for $(g(t),f(t))$, which is denoted by $S_n(x)\sim (g(t),f(t))$. By \eqref{8}, we easily get that $\left. \left<e^{yt}\right| p(x)\right>=p(y)$. For $f(t)\in{\mathcal{F}}$ and $p(x)\in{\mathbb{P}}$, we have
\begin{equation}\label{9}
f(t)=\sum_{k=0} ^{\infty} \frac{\left<f(t)|x^k\right>}{k!} t^k,~p(x)=\sum_{k=0} ^{\infty} \frac{\left<t^k|p(x)\right>}{k!}x^k,
\end{equation}
and
\begin{equation}\label{10}
\left. \left<f_1(t)\cdots f_m(t)\right|x^n \right>=\sum_{i_1+\cdots+i_m=n}\binom{n}{i_1,\ldots,i_m}\left(\prod_{j=1} ^m \left.\left<f_j(t)\right|x^{i_j}\right>\right),
\end{equation}
where $f_1(t),f_2(t),\ldots,f_m(t) \in {\mathcal{F}}$ (see \cite{09,10,18,19}).
For $f(t),g(t)\in {\mathcal{F}}$ and $p(x)\in{\mathbb{P}}$, by \eqref{9}, we get
\begin{equation}\label{11}
p^{(k)}(0)=\left<t^k|p(x)\right>,~\left<1\left|p^{(k)}(x)\right.\right>=p^{(k)}(0).
\end{equation}
Thus, by \eqref{11}, we have
\begin{equation}\label{12}
t^kp(x)=p^{(k)}(x)=\frac{d^kp(x)}{dx^k},~(k \geq 0), {\text{ (see [10,11,18,19])}}.
\end{equation}

Let $S_n(x) \sim \left(g(t),f(t)\right)$. Then we have
\begin{equation}\label{12-1}
\frac{1}{g({\bar{f}}(t))}e^{y{\bar{f}}(t)}=\sum_{k=0} ^{\infty} \frac{S_k(y)}{k!}t^k,{\text{ for all }}y \in {\mathbb{C}},
\end{equation}
where ${\bar{f}}(t)$ is the compositional inverse of $f(t)$ (see \cite{18,19}). By \eqref{2} and \eqref{12-1}, we see that $N_n ^{(a)}(x)\sim\left(\left(\frac{e^t-1}{t}\right)^a,e^t-1\right)$.

For $a \neq 0$, the Poisson-Charlier sequences are given by
\begin{equation}\label{13}
C_n(x;a)=\sum_{k=0} ^n \binom{n}{k}(-1)^{n-k} a^{-k}(x)_k\sim \left(e^{a(e^t-1)},a(e^t-1)\right).
\end{equation}
In particular, $n\in{\mathbb{Z}}_+={\mathbb{N}}\cup\left\{0\right\}$, we have
\begin{equation}\label{14}
\sum_{l=0} ^{\infty} C_n(l;a)\frac{t^l}{l!}=e^t\left(\frac{t-a}{a}\right)^n,{\text{ (see [18,19])}}.
\end{equation}
The {\it{Frobenius-type Eulerian polynomials}} of order $a$ are given by
\begin{equation}\label{15}
\left(\frac{1-\lambda}{e^{t(\lambda-1)}-\lambda}\right)^ae^{xt}=\sum_{n=0} ^{\infty}A_n ^{(a)} (x|\lambda),{\text{ (see [17,18])}}.
\end{equation}
From \eqref{12-1} and \eqref{15}, we note that
\begin{equation*}
A_n ^{(a)} (x|\lambda)\sim\left(\left(\frac{e^{t(1-\lambda)}-\lambda}{1-\lambda}\right)^a,t\right).
\end{equation*}
Let us assume that $p_n(x)\sim (1,f(t))$, $q_n(x)\sim(1,g(t))$. Then we have
\begin{equation}\label{16}
q_n(x)=x\left(\frac{f(t)}{g(t)}\right)^nx^{-1}p_n(x),{\text{ (see [18,19])}}.
\end{equation}
The equation \eqref{16} is important in deriving our results in this paper. The purpose of this paper is to investigate some properties of Sheffer sequences of several polynomials arising from umbral calculus. From our investigation, we can derive many interesting identities of several polynomials.

\section{Sheffer sequences of polynomials}

Let us assume that $S_n(x)\sim(g(t),f(t))$. Then, by the definition of Sheffer sequence, we see that $g(t)S_n(x)\sim(1,f(t))$. If $g(t)$ is an invertible series, then $\frac{1}{g(t)}$ is also an invertible series. Let us consider the following Sheffer sequences:
\begin{equation}\label{17}
M_n(x)\sim(1,f(t)),~x^n\sim(1,t).
\end{equation}
From \eqref{16} and \eqref{17}, we note that
\begin{equation}\label{18}
M_n(x)=x \left(\frac{t}{f(t)}\right)^nx^{-1}x^n=x\left(\frac{t}{f(t)}\right)^nx^{n-1}.
\end{equation}
For $g(t)S_n(x)\sim(1,f(t))$, by \eqref{18}, we get
\begin{equation}\label{19}
g(t)S_n(x)=x\left(\frac{t}{f(t)}\right)^nx^{n-1}.
\end{equation}
Therefore, by \eqref{19}, we obtain the following theorem.
\begin{theorem}\label{thm1}
For $S_n(x)\sim(g(t),f(t))$ and $n \geq 1$, we have
\begin{equation*}
S_n(x)=\frac{1}{g(t)}x\left(\frac{t}{f(t)}\right)^nx^{n-1}.
\end{equation*}
\end{theorem}
For example, let $S_n(x)=D_n(x)\sim\left(\frac{1-\lambda}{e^t-\lambda},\frac{e^t-1}{e^t+1}\right)$, where $D_n(x)$ is the $n$-th Daehee polynomial (see \cite{01,08,09}). Then, by Theorem \ref{thm1}, we get
\begin{equation*}
\begin{split}
D_n(x)&=\left(\frac{e^t-\lambda}{1-\lambda}\right)x\left(\frac{t}{e^t-1}\right)^n(e^t+1)^nx^{n-1}=\left(\frac{e^t-\lambda}{1-\lambda}\right)x \sum_{l=0} ^n \binom{n}{l}B_{n-1} ^{(n)} (x+l) \\
&=\frac{1}{1-\lambda}\sum_{l=0} ^n \binom{n}{l}\left\{(x+1)B_{n-1} ^{(n)} (x+l+1)-\lambda x B_{n-1} ^{(n)} (x+l)\right\}.
\end{split}
\end{equation*}

Let us take $S_n(x)\sim\left(\left(\frac{e^t-\lambda}{1-\lambda}\right)^a,\frac{t^2}{e^{bt}-1}\right)$, $(b\neq0)$. Then, by Theorem \ref{thm1}, we get
\begin{equation}\label{20}
\begin{split}
S_n(x)&=\left(\frac{1-\lambda}{e^t-\lambda}\right)^ax\left(\frac{e^{bt}-1}{t}\right)^nx^{n-1} \\
&=\left(\frac{1-\lambda}{e^t-\lambda}\right)^ax\sum_{k=0} ^{n-1} \frac{n!b^{k+n}}{(k+n)!} S_2(k+n,n)x^{n-k-1}(n-1)_k \\
&=\sum_{k=0} ^{n-1} \frac{\binom{n-1}{k}}{\binom{k+n}{n}}S_2(k+n,n)b^{k+n}H_{n-k} ^{(a)}(x|\lambda).
\end{split}
\end{equation}
Therefore, by \eqref{20}, we obtain the following theorem.
\begin{theorem}\label{thm2}
For $n \geq 1$, let $S_n(x)\sim\left(\left(\frac{e^t-\lambda}{1-\lambda}\right)^a,\frac{t^2}{e^{bt}-1}\right),~b\neq0$. Then we have
\begin{equation*}
S_n(x)=\sum_{k=0} ^{n-1} \frac{\binom{n-1}{k}}{\binom{k+n}{n}}S_2(k+n,n)b^{k+n}H_{n-k} ^{(a)}(x|\lambda).
\end{equation*}
\end{theorem}
Let
\begin{equation}\label{21}
S_n(x)\sim\left(\left(\frac{e^t-1}{t}\right)^a,\frac{t^2e^{bt}}{e^{ct}-1}\right),~c\neq0.
\end{equation}
From Theorem \ref{thm1}, we can derive
\begin{equation}\label{22}
\begin{split}
& S_n(x) \\
=&\left(\frac{t}{e^t-1}\right)^ax\left(\frac{e^{ct}-1}{te^{bt}}\right)^nx^{n-1} \\
=&\left(\frac{t}{e^t-1}\right)^axe^{-nbt}\sum_{l=0} ^{\infty}\frac{n!S_2(l+n,n)}{(l+n)!}c^{l+n}t^lx^{n-1}\\
=&\left(\frac{t}{e^t-1}\right)^ax\sum_{l=0} ^{n-1} \frac{\binom{n-1}{l}}{\binom{l+n}{l}}S_2(l+n,n)c^{n+l}(x-nb)^{n-1-l}\\
=&\left(\frac{t}{e^t-1}\right)^ax\sum_{l=0} ^{n-1} \sum_{j=0} ^{n-1-l}\frac{\binom{n-1}{l}}{\binom{l+n}{l}}\binom{n-1-l}{j}S_2(l+n,n)c^{n+l}(-nb)^jx^{n-1-l-j}\\
=&\sum_{l=0} ^{n-1} \sum_{j=0} ^{n-1-l}\frac{\binom{n-1}{l}}{\binom{l+n}{l}}\binom{n-1-l}{j}S_2(l+n,n)c^{n+l}(-nb)^jB_{n-l-j} ^{(a)} (x).
\end{split}
\end{equation}
Therefore, by \eqref{22}, we obtain the following theorem.
\begin{theorem}\label{thm3}
For $n \geq 1$, let $S_n(x)\sim\left(\left(\frac{e^t-1}{t}\right)^a,\frac{t^2e^{bt}}{e^{ct}-1}\right),~c\neq0$. Then we have
\begin{equation*}
S_n(x)=\sum_{l=0} ^{n-1} \sum_{j=0} ^{n-1-l}\frac{\binom{n-1}{l}}{\binom{l+n}{l}}\binom{n-1-l}{j}S_2(l+n,n)c^{n+l}(-nb)^jB_{n-l-j} ^{(a)} (x).
\end{equation*}
\end{theorem}
Let us take the following Sheffer sequence:
\begin{equation}\label{23}
S_n(x)\sim\left(\left(\frac{e^t+1}{2}\right)^{\alpha},\frac{t^2}{\log (1+t)}\right).
\end{equation}
By Theorem \ref{thm1} and \eqref{23}, we get
\begin{equation}\label{24}
\begin{split}
& S_n(x) \\
=&\left(\frac{2}{e^t+1}\right)^{\alpha}x\left(\frac{\log (1+t)}{t}\right)^nx^{n-1} =\left(\frac{2}{e^t+1}\right)^{\alpha}x\sum_{l=0} ^{\infty} \frac{N_l ^{(n)}}{l!}t^lx^{n-1} \\
=&\left(\frac{2}{e^t+1}\right)^{\alpha}x\sum_{l=0} ^{n-1} \binom{n-1}{l} N_l ^{(n)} x^{n-l-1} \\
=&\sum_{l=0} ^{n-1} \binom{n-1}{l} N_l ^{(n)} E_{n-l} ^{(\alpha)} (x),
\end{split}
\end{equation}
where $E_n ^{(\alpha)} (x)$ are the $n$-th Euler polynomials of order $\alpha$ which is defined by the generating function to be
\begin{equation*}
\left(\frac{2}{e^t+1}\right)^{\alpha}e^{xt}=\sum_{n=0} ^{\infty} E_n ^{(\alpha)} (x) \frac{t^n}{n!}.
\end{equation*}
Therefore, by \eqref{24}, we obtain the following theorem.
\begin{theorem}\label{thm4}
For $n \geq 1$, let $S_n(x)\sim\left(\left(\frac{e^t+1}{2}\right)^{\alpha},\frac{t^2}{\log (1+t)}\right)$.
Then we have
\begin{equation*}
S_n(x)=\sum_{l=0} ^{n-1} \binom{n-1}{l} N_l ^{(n)} E_{n-l} ^{(\alpha)} (x).
\end{equation*}
\end{theorem}
As is known, we note that
\begin{equation}\label{25}
\left(\frac{\log (1+t)}{t}\right)^n=n\sum_{l=0} ^{\infty} \frac{B_l ^{(n+l)}}{n+l} \frac{t^l}{l!}.
\end{equation}
Thus, by Theorem \ref{thm1} and \eqref{35}, we get
\begin{equation}\label{26}
\begin{split}
& S_n(x) \\
=&\left(\frac{2}{e^t+1}\right)^{\alpha}x\left(\frac{\log (1+t)}{t}\right)^nx^{n-1} \\
=&\left(\frac{2}{e^t+1}\right)^{\alpha}xn\sum_{l=0} ^{n-1} \frac{B_l ^{(n+l)}}{n+l}\binom{n-1}{l} x^{n-1-l} \\
=&n\sum_{l=0} ^{n-1} \frac{B_l ^{(n+l)}}{n+l}\binom{n-1}{l} E_{n-l} ^{(\alpha)} (x).
\end{split}
\end{equation}
Therefore, by Theorem \ref{thm4} and \eqref{26}, we obtain the following corollary.
\begin{corollary}\label{coro5}
For $n \geq 1$, and $0 \leq l \leq n-1$, we have
\begin{equation*}
\frac{N_l ^{(n)}}{n}=\frac{B_l ^{(n+l)}}{n+l}.
\end{equation*}
\end{corollary}

{\scshape Remark.} Let $S_n(x)\sim\left(\left(\frac{e^t-1}{t}\right)^{\alpha},\log (1+t)\right)$. Then, by Theorem \ref{thm1}, we get
\begin{equation}\label{27}
\begin{split}
S_n(x)&=\left(\frac{t}{e^t-1}\right)^{\alpha}x\left(\frac{t}{\log (1+t)}\right)^n x^{n-1}\\
&=\left(\frac{t}{e^t-1}\right)^{\alpha}x\sum_{l=0} ^{n-1} \binom{n-1}{l} N_l ^{(-n)} x^{n-1-l} \\
&=\sum_{l=0} ^{n-1} \binom{n-1}{l} N_l ^{(-n)} B_{n-l} ^{(\alpha)} (x).
\end{split}
\end{equation}
Let us assume that
\begin{equation}\label{28}
S_n(x)\sim\left(\left(\frac{e^t-\lambda}{1-\lambda}\right)^{\alpha},\frac{\log (1+t)}{(1+t)^c}\right),~(c \neq 0).
\end{equation}
Then, by Theorem \ref{thm1} and \eqref{28}, we get
\begin{equation}\label{29}
\begin{split}
& S_n(x) \\
=&\left(\frac{1-\lambda}{e^t-\lambda}\right)^{\alpha}x\left(\frac{t(1+t)^c}{\log (1+t)}\right)^nx^{n-1} \\
=&\left(\frac{1-\lambda}{e^t-\lambda}\right)^{\alpha}x\sum_{l=0} ^{n-1} B_l ^{(l-n+1)} (cn+1)\frac{(n-1)_l}{l!}x^{n-1-l} \\
=&\sum_{l=0} ^{n-1}\binom{n-1}{l}B_l ^{(l-n+1)}(cn+1)\left(\frac{1-\lambda}{e^t-\lambda}\right)^{\alpha}x^{n-l} \\
=&\sum_{l=0} ^{n-1} \binom{n-1}{l}B_l ^{(l-n+1)}(cn+1) H_{n-l} ^{(\alpha)} (x|\lambda).
\end{split}
\end{equation}
Therefore, by \eqref{29}, we obtain the following theorem.
\begin{theorem}\label{thm6}
For $n \geq 1$, let $S_n(x)\sim\left(\left(\frac{e^t-\lambda}{1-\lambda}\right)^{\alpha},\frac{\log (1+t)}{(1+t)^c}\right)$, $c\neq 0$. Then we have
\begin{equation*}
S_n(x)=\sum_{l=0} ^{n-1} \binom{n-1}{l}B_l ^{(l-n+1)}(cn+1) H_{n-l} ^{(\alpha)} (x|\lambda).
\end{equation*}
\end{theorem}

As is well known, the {\it{Bernoulli polynomials of the second kind}} are defined by the generating function to be
\begin{equation}\label{30}
\frac{t(1+t)^x}{\log (1+t)}=\sum_{l=0} ^{\infty} \frac{b_l (x)}{l!}t^l,{\text{ (see [18,19])}}.
\end{equation}
Thus, by \eqref{10} and \eqref{30}, we get
\begin{equation}\label{31}
\left(\frac{t(1+t)^c}{\log (1+t)}\right)^n=\sum_{l=0} ^{\infty}\left(\sum_{l_1+\cdots+l_n=l}\binom{l}{l_1,\ldots,l_n}b_{l_1} (c)\cdots b_{l_n} (c)\right)\frac{t^l}{l!}.
\end{equation}
By Theorem \ref{thm1},\eqref{28} and \eqref{31}, we get
\begin{equation}\label{32}
\begin{split}
&S_n(x) \\
=&\left(\frac{1-\lambda}{e^t-\lambda}\right)^{\alpha}x\sum_{l=0} ^{n-1} \left(\sum_{l_1+\cdots+l_n=l} \binom{l}{l_1,\ldots,l_n}\left(\prod_{i=1} ^n b_{l_i} (c) \right)\binom{n-1}{l}x^{n-1-l}\right) \\
=&\sum_{l=0} ^{n-1} \left(\sum_{l_1+\cdots+l_n=l}\binom{l}{l_1,\ldots,l_n}\left(\prod_{i=1} ^n b_{l_i} (c) \right)\right) \binom{n-1}{l} \left(\frac{1-\lambda}{e^t-\lambda}\right)^{\alpha}x^{n-l} \\
=&\sum_{l=0} ^{n-1} \left(\sum_{l_1+\cdots+l_n=l}\binom{l}{l_1,\ldots,l_n}\left(\prod_{i=1} ^n b_{l_i} (c) \right)\right)\binom{n-1}{l}H_{n-l} ^{(\alpha)}(x|\lambda).
\end{split}
\end{equation}
Therefore, by Theorem \ref{thm6} and \eqref{32}, we obtain the following theorem.
\begin{theorem}\label{thm7}
For $n\geq 1$, $0 \leq l \leq n-1$, we have
\begin{equation*}
\sum_{l_1+\cdots+l_n=l}\binom{l}{l_1,\ldots,l_n}\left(\prod_{i=1} ^n b_{l_i} (c) \right)=B_l ^{(l-n+1)}(cn+1),~(c\neq 0).
\end{equation*}
\end{theorem}

{\scshape Remark.} From \eqref{2}, we note that
\begin{equation}\label{33}
\left(\frac{t(1+t)^c}{\log (1+t)}\right)^n x^{n-1}=\sum_{l=0} ^{n-1} \binom{n-1}{l}N_l ^{(-n)}(cn)x^{n-1-l},
\end{equation}
where $c \neq 0$. By Theorem \ref{thm1}, \eqref{28} and \eqref{33}, we get
\begin{equation}\label{34}
\begin{split}
S_n(x)=&\left(\frac{1-\lambda}{e^t-\lambda}\right)^{\alpha}x\left(\frac{t(1+t)^c}{\log (1+t)}\right)^nx^{n-1} \\
&=\sum_{l=0} ^{n-1} \binom{n-1}{l}N_l ^{(-n)} (cn)H_{n-l} ^{(\alpha)}(x|\lambda).
\end{split}
\end{equation}
From \eqref{32} and \eqref{34}, we can derive the following identity:
\begin{equation}\label{35}
N_l ^{(-n)}(cn)=\sum_{l_1+\cdots+l_n=l}\binom{l}{l_1,\ldots,l_n}\left(\prod_{i=1} ^n b_{l_i} (c)\right),
\end{equation}
where $n\geq1$, $0 \leq l \leq n-1$ and $c \neq 0$. Let
\begin{equation}\label{36}
S_n(x)\sim\left(\left(\frac{e^{(\lambda-1)t}-\lambda}{1-\lambda}\right)^{\alpha},\frac{t^2(1+t)^c}{\log (1+t)}\right),~c \neq 0.
\end{equation}
From Theorem \ref{thm1} and \eqref{36}, we note that
\begin{equation}\label{37}
\begin{split}
S_n(x)=&\left(\frac{1-\lambda}{e^{(\lambda-1)t}-\lambda}\right)^{\alpha}x\left(\frac{\log (1+t)}{t(1+t)^c}\right)^nx^{n-1} \\
=& \left(\frac{1-\lambda}{e^{(\lambda-1)t}-\lambda}\right)^{\alpha}x\sum_{l=0} ^{n-1} \binom{n-1}{l} N_l ^{(n)} (-cn) x^{n-1-l} \\
=&\sum_{l=0} ^{n-1} \binom{n-1}{l} N_l ^{(n)} (-cn)A_{n-l} ^{(\alpha)} (x|\lambda).
\end{split}
\end{equation}
Therefore, by \eqref{37}, we obtain the following proposition.
\begin{proposition}\label{prop8}
For $n \geq 1$, let $S_n(x)\sim \left(\left(\frac{e^{(\lambda-1)t}-\lambda}{1-\lambda}\right)^{\alpha},\frac{t^2(1+t)^c}{\log (1+t)}\right)$, $c \neq 0$. Then we have
\begin{equation*}
S_n(x)=\sum_{l=0} ^{n-1} \binom{n-1}{l}N_l ^{(n)} (-nc) A_{n-l} ^{(\alpha)} (x|\lambda) .
\end{equation*}
\end{proposition}
Now we observe that
\begin{equation}\label{38}
\begin{split}
\left(\frac{\log (1+t) }{t(1+t)^c}\right)^n=&(1+t)^{-nc}\left(\frac{\log (1+t)}{t}\right)^n \\
=&(1+t)^{-nc}\left(\sum_{k=0} ^{\infty} \frac{n!S_1(k+n,n)}{(k+n)!}t^k\right) \\
=&\left(\sum_{m=0} ^{\infty} \binom{-nc}{m}t^m\right)\left(\sum_{k=0} ^{\infty} \frac{n!S_1(k+n,n)}{(k+n)!}t^k\right) \\
=&\sum_{l=0} ^{\infty}\left\{\sum_{k=0} ^l \frac{n!S_1(k+n,n)}{(k+n)!}\binom{-nc}{l-k}\right\}t^l.
\end{split}
\end{equation}
By Theorem \ref{thm1}, \eqref{36} and \eqref{38}, we get
\begin{equation}\label{39}
\begin{split}
S_n(x)=&\left(\frac{1-\lambda}{e^{(\lambda-1)t}-\lambda}\right)^{\alpha}x\left(\frac{\log (1+t)}{t(1+t)^c}\right)^n x^{n-1} \\
=&\sum_{l=0} ^{n-1} \binom{n-1}{l}l!\left\{\sum_{k=0} ^l \frac{n!}{(k+n)!}S_1(n+k,n)\binom{-nc}{l-k}\right\}A_{n-l} ^{(\alpha)} (x|\lambda).
\end{split}
\end{equation}
Therefore, by Proposition \ref{prop8} and \eqref{39}, we obtain the following theorem.
\begin{theorem}\label{thm9}
For $n \geq 1$, $0 \leq l \leq n-1$ and $c \neq 0$, we have
\begin{equation*}
N_l ^{(n)}(-cn)=l!\sum_{k=0} ^l \frac{n!}{(n+k)!}S_1(k+n,n)\binom{-nc}{l-k}.
\end{equation*}
\end{theorem}

{\scshape Remark.} It is easy to show that
\begin{equation}\label{40}
\left(\log (1+t)\right)^n=\sum_{l=0} ^{\infty} \frac{n!}{(l+n)!}S_1(l+n,k)t^{l+n}.
\end{equation}
By Theorem \ref{thm1}, \eqref{23} and \eqref{40}, we get
\begin{equation}\label{41}
\begin{split}
S_n(x)=&\left(\frac{2}{e^t+1}\right)^{\alpha}x\left(\frac{\log (1+t)}{t}\right)^nx^{n-1} \\
=&\left(\frac{2}{e^t+1}\right)^{\alpha}x\sum_{l=0} ^{n-1} \frac{n!l!}{(l+n)!}\binom{n-1}{l}S_1(l+n,n)x^{n-1-l}\\
=&\sum_{l=0} ^{n-1} \frac{\binom{n-1}{l}}{\binom{l+n}{n}}S_2(l+n,n)E_{n-l} ^{(\alpha)}(x).
\end{split}
\end{equation}
From Theorem \ref{thm4} and \eqref{41}, we can derive the following identity:
\begin{equation}\label{42}
N_l ^{(n)}=\frac{S_2(l+n,n)}{\binom{l+n}{n}},{\text{ where }}n \geq 1,~0 \leq l \leq n-1.
\end{equation}
Let us consider the following Sheffer sequence:
\begin{equation}\label{43}
S_n(x)\sim\left(\left(\frac{e^{(\lambda-1)t}-\lambda}{1-\lambda}\right)^{\alpha},\frac{t}{e^{ct}(1+bt)^m}\right),~b,c\neq 0,~m \in {\mathbb{Z}}_+.
\end{equation}
By Theorem \ref{thm1} and \eqref{43}, we get
\begin{equation}\label{44}
\begin{split}
S_n(x)=&\left(\frac{1-\lambda}{e^{(\lambda-1)t}-\lambda}\right)^{\alpha}x\left(e^{ct}(1+bt)^m\right)^nx^{n-1}\\
=&\left(\frac{1-\lambda}{e^{(\lambda-1)t}-\lambda}\right)^{\alpha}xe^{nct}(1+bt)^{mn}x^{n-1}.
\end{split}
\end{equation}
From \eqref{14} and \eqref{44}, we can derive
\begin{equation}\label{45}
\begin{split}
S_n(x)=&\left(\frac{1-\lambda}{e^{(\lambda-1)t}-\lambda}\right)^{\alpha}x(-1)^{mn}\sum_{l=0} ^{n-1}C_{mn}\left(l;-\frac{nc}{b}\right)(nc)^l\binom{n-1}{l}x^{n-1-l} \\
=&(-1)^{mn}\sum_{l=0} ^{n-1}C_{mn}\left(l;-\frac{nc}{b}\right)(nc)^l\binom{n-1}{l}A_{n-l} ^{(\alpha)}(x|\lambda).
\end{split}
\end{equation}
Therefore, by \eqref{45}, we obtain the following theorem.
\begin{theorem}
For $n \geq 1$, let $S_n(x)\sim\left(\left(\frac{e^{(\lambda-1)t}-\lambda}{1-\lambda}\right)^{\alpha},\frac{t}{e^{ct}(1+bt)^m}\right)$, where $m \in {\mathbb{Z}}_+$, $b \neq 0$ and $c\neq 0$. Then we have
\begin{equation*}
S_n(x)=(-1)^{mn}\sum_{l=0} ^{n-1}C_{mn}\left(l;-\frac{nc}{b}\right)(nc)^l\binom{n-1}{l}A_{n-l} ^{(\alpha)}(x|\lambda).
\end{equation*}
\end{theorem}


\begin{thebibliography}{10}

\bibitem {01} L. Carlitz, {\it Eulerian numbers and polynomials of higher order}, Duke Math. J.,  ${\mathbf{27}}$  (1960),  401-423.

\bibitem {02} B. Diarra, {\it Ultrametric umbral calculus in characteristic $p$}, Bull. Belg. Math. Soc. Simon Stevin,  ${\mathbf{14}}$ (2007) 845-869.

\bibitem {03} R. Dere, Y. Simsek, {\it Applications of umbral algebra to some special polynomials},  Adv. Stud. Contemp. Math. ${\mathbf{22}}$ (2012) 433-438.

\bibitem {04} T. Ernst, {\it Examples of a $q$-umbral calculus}, Adv. Stud. Contemp. Math. ${\mathbf{16}}$ (2008), no. 1, 1-22.


\bibitem {05} T. Kim,  {\it Some identities on the $q$-Euler polynomials of higher order and $q$-Stirling numbers by the fermionic $p$-adic integral on ${\mathbb{Z}}_p$}, Russ. J. Math. Phys., ${\mathbf{16}}$   (2009),  no. 4, 484-491.

\bibitem {06} T. Kim,  {\it Identities involving Frobenius-Euler polynomials arising from non-linear differential equations}, J. Number Theory, ${\mathbf{132}}$   (2012),  no. 1, 2854-2865.

\bibitem {07} T. Kim,  {\it An identity of the symmetry for the Frobenius-Euler polynomials associated with the fermionic $p$-adic invariant $q$-integrals on ${\mathbb{Z}}_p$}, Rocky Mountain J. Math., ${\mathbf{41}}$   (2011),  no. 1, 239-247.

\bibitem {08} T. Kim,  {\it Symmetry $p$-adic invariant integral on ${\mathbb{Z}}_p$ for Bernoulli and Euler polynomials}, J. Difference Equ. Appl., ${\mathbf{14}}$   (2008),  no. 12, 1267-1277.

\bibitem {09} D. S. Kim, T. Kim, S. H. Lee and S. H. Rim, {\it Frobenius-Euler polynomials and umbral calculus in the $p$-adic case}, Adv. Difference Equ. 2012, 2012:222.

\bibitem {10} D. S. Kim and T. Kim, {\it Some new identities of Frobenius-Euler numbers and polynomials},  J. of Inequ. and Appl., 2012, 2012:307.

\bibitem {11} D. S. Kim and T. Kim, {\it Applications of Umbral Calculus Associated with $p$-Adic Invariant Integrals on ${\mathbb{Z}}_p$}, Abstract and Applied Analysis 2012 (2012), Article ID 865721, 12 pages.

\bibitem {12} D. S. Kim and T. Kim, {\it Some identities of Frobenius-Euler polynomials arising from umbral calculus}, Adv. Difference Equ., 2012, 2012:196.

\bibitem {13} D. S. Kim, T. Kim, S.-H. Lee and Y.-H. Kim, {\it Some identities for the product of two Bernoulli and Euler polynomials}, Adv. Difference Equ., 2012, 2012:95, 14 pp.

\bibitem {14} T. Kim, S.-H. Rim, D. V. Dolgy and S.-H. Lee, {\it Some identities on Bernoulli and Euler polynomials arising from the orthogonality of Laguerre polynomials}, Adv. Difference Equ., 2012, 2012:201.

\bibitem {15} T. Kim, {\it Some identities on the $q$-Euler polynomials of higher order and $q$-Stirling numbers by the fermionic $p$-adic integral on ${\mathbb{Z}}_+$}, Russ. J. Math. Phys., ${\mathbf{16}}$ (2009), no. 4, 484-491.

\bibitem {16} T. Mansour, M. Schork, S. Severini, {\it A generalization of boson normal ordering}, Phys. Lett. A, ${\mathbf{364}}$ (2007), no. 3-4, 214-220.

\bibitem {17}  T. J. Robinson, {\it Formal calculus and umbral calculus},  Electron. J. Combin., ${\mathbf{17}}$ (2010), Research Paper 95, 31 pp.

\bibitem {18} S. Roman, {\it More on the umbral calculus, with emphasis on the $q$-umbral calculus}, J. Math. Anal. Appl., ${\mathbf{107}}$ (1985), 222-254.

\bibitem {19} S. Roman, {\it The umbral calculus}, Dover Publ. Inc. New York, 2005.

\bibitem {20} C. Ryoo, {\it Some relations between twisted q-Euler numbers and Bernstein polynomials}, Adv. Stud. Contemp. Math. ${\mathbf{21}}$ (2011), no. 2, 217-223.
\bibitem {21} S. Araci, M. Acikgoz, {\it  A note on the Frobenius-Euler numbers and polynomials associated with Bernstein polynomials},
Adv. Stud. Contemp. Math., ${\mathbf{22}}$ (2012), no. 3, 399-406.

\end{thebibliography}
\end{document}